\documentclass[10pt,twoside]{article}
\usepackage{amsmath,amsthm,amsfonts,latexsym,amscd}

\usepackage{amsthm,amssymb}

\swapnumbers
\theoremstyle{plain}
\newtheorem{theoreme}{Th\'{e}or\`{e}me}
\newtheorem{theorem}{Theorem}
\newtheorem{lemma}[theorem]{Lemma}
\newtheorem{corollary}[theorem]{Corollary}
\newtheorem{proposition}[theorem]{Proposition}

\newtheorem{too_optimistic_conjecture}[theorem]{The Strong Atiyah Conjecture}
\theoremstyle{definition}

\theoremstyle{remark}
\newtheorem{remark}[theorem]{Remark}

\usepackage{amsmath,amsfonts,latexsym}

\newcommand{\reals}{\mathbb{R}}

\newcommand{\naturals}{\mathbb{N}}

\newcommand{\integers}{\mathbb{Z}}

\newcommand{\rationals}{\mathbb{Q}}

\newcommand{\abs}[1]{\left\lvert#1\right\rvert} 

\newcommand{\generate}[1]{\langle#1\rangle}
\newcommand{\tensor}{\otimes}
\newcommand{\into}{\hookrightarrow}
\newcommand{\onto}{\twoheadrightarrow}
\newcommand{\iso}{\cong}

\newcommand{\semiProd}{\rtimes}

\DeclareMathOperator{\tr}{tr}
\DeclareMathOperator{\pr}{pr}

\newcommand{\forget}[1]{}

{\catcode`@=11\global\let\c@equation=\c@theorem}

\allowdisplaybreaks[2]

\DeclareMathOperator{\fin}{fin^{-1}} 

\usepackage{amssymb,latexsym,epic,eepic,epsfig,psfig}

\title{On a conjecture of Atiyah }

\author{
Rostislav I. GRIGORCHUK\\  
{\small Steklov Mathematical Institute, Gubkina Str.~8 Moscow, 117966, Russia}\\
{\small e-mail: grigorch@mi.ras.ru }\\ 
Peter LINNELL\\
{\small Department of Mathematics,
              Virginia Polytechnic Institute and State University}\\
{\small            Blacksburg, VA 24061, USA}\\
{\small e-mail: linnell@math.vt.edu}\\     
Thomas SCHICK\\
{\small FB Mathematik, Universit{\"a}t M{\"u}nster, 
Einsteinstr.~62, 48159 M{\"unster}, Germany}\\
{\small e-mail: thomas.schick@math.uni-muenster.de
}\\
Andrzej \.ZUK\\
{\small   CNRS, Ecole Normale Sup\'erieure de Lyon, Unit\'e de  Math\'ematiques Pures et Appliqu\'ees}\\
{\small 46, All\'ee d'Italie,  F-69364 Lyon cedex 07,  France}\\
{\small e-mail: azuk@umpa.ens-lyon.fr }
}
\date{}

\begin{document}
\maketitle

\noindent {\small {\bf Abstract.} 
In this note we explain how the computation of the spectrum of the lamplighter
group from \cite{Grigorchuk-Zuk(2000)} yields a counterexample to a strong
version of the Atiyah conjectures 
about the range of $L^2$-Betti numbers of closed manifolds. 

\smallskip
\noindent MSC-number: 57N65}

\bigskip

\centerline{\large \bf Sur une conjecture d'Atiyah}

\bigskip

\noindent {\small {\bf R\'esum\'e.}
Dans cette note on montre comment le calcul du 
spectre du groupe de l'allumeur de
r\'everb\`eres fait dans \cite{Grigorchuk-Zuk(2000)} donne un contre-exemple 
\`a une des conjectures d'Atiyah
sur les nombres de Betti $L^2$ des vari\'et\'es ferm\'ees.
}

\medskip

{\large{\bf Version fran{\c c}aise abr\'eg\'ee}}

\medskip

Dans \cite{Atiyah(1976)} Atiyah a demand\'e si les nombres de Betti $L^2$ 
des rev{\^e}tements universels des vari\'et\'es ferm\'ees sont 
toujours rationnels. Plus tard cette question 
a \'et\'e transform\'ee en conjecture  sous le nom d'Atiyah 
(\cite[Theorem 8]{Cohen(1979)}, \cite[Conjecture 7.1]{Lott-Lueck(1995)}, 
\cite[Conjecture
9.8]{Linnell(1998)},
\cite[Conjecture 5.1]{Lueck(1998)},
\cite[Conjecture 5.1]{Reich(1999)},
\cite[Conjecture 13]{Lueck(2000)}, \cite[Definition
1.2]{Dicks-Schick(2000)}).
La r\'eponse ne d\'epend que du groupe
 fondamental.  
La plus c\'el\`ebre parmi ces conjectures
est la suivante~: les nombres de Betti $L^2$ d'un groupe de pr\'esentation
finie sans torsion sont entiers. La conjecture  d'Atiyah forte dit
que les nombres de Betti $L^2$ d'un groupe de pr\'esentation finie $G$
sont des rationnels dont les d\'enominateurs 
sont d\'etermin\'es par les ordres des sous-groupes
finis de $G$.
Il y a aussi des versions encore plus fortes de la conjecture d'Atiyah, quand
le groupe n'est pas suppos\'e  \^etre de pr\'esentation finie 
\cite{Lueck(1997)}.
 Il y 
a des reformulation equivalente de cette conjecture en terme de $CW$
complexes et en terme des vari\'et\'es.

Plusieurs r\'esultats confirment diff\'erentes formes de la 
conjecture d'Atiyah  (\cite{Cohen(1979)}, \cite{Linnell(1993)},
\cite{Schick(1999b)}, \cite{Dicks-Schick(2000)}). 
Un de ces resultats est le th\'eor\`eme de Linnell
\cite{Linnell(1993)} qui affirme que la conjecture d'Atiyah forte est 
vraie pour les groupes 
moyennables \'el\'ementaires avec une borne uniforme sur
l'ordre des sous-goupes finis et m\^eme pour une famille plus
large de groupes contenant les groupes libres. D'autre part 
les resultats de \cite{Dicks-Schick(2000)} essentiallement montrent que la famille des groupes
pour lesquels la conjecture  d'Atiyah forte est vraie est ferm\'ee
sous les HNN extensions et les produits  amalgam\'es
sous l'hypoth\`ese que les ordres des sous-groupes finis de ces groupes
sont uniform\'ement born\'es.

Le r\'esultat principal de cet article montre que la conjecture
d'Atiyah forte est fausse.

\begin{theoreme}
Soit $G$ le groupe donn\'e par la pr\'esentation
$$G = \langle a,t,s \mid a^2 =1 , [t,s]=1 , 
[t^{-1}at,a] = 1, s^{-1}as= at^{-1}at \rangle.$$
Le groupe $G$ est 
metab\'elien et donc moyennable \'el\'ementaire.
Chaque sous-groupe fini de $G$ est un 2-groupe abelien \'el\'ementaire,
en particulier l'ordre de chaque sous-groupe fini de $G$ est 
une puissance de 2.
Il existe une vari\'et\'e riemannienne ferm\'ee $(M,g)$ 
de dimension 7 telle que $\pi_1(M) =G$ pour laquelle le troisi\`eme
nombre de Betti $L^2$ est \'egal \`a
$$ b_3^{(2)}(M,g) =\frac{1}{3}.$$
\end{theoreme}

Ce th\'eor\`eme montre qu'on ne peut pas g\'en\'eraliser les r\'esultats 
de Linnell \cite{Linnell(1993)} et Dick-Schick \cite{Dicks-Schick(2000)}
aux groupes moyennables \'el\'ementaires sans 
une borne uniforme sur l'ordre des sous-groupes finis. 

La preuve du Th\'eor\`eme 1 est bas\'ee sur les r\'esultats de 
\cite{Grigorchuk-Zuk(2000)} sur le
spectre et la mesure spectrale de l'op\'erateur de Markov $A$ de la marche 
al\'eatoire simple sur le groupe de l'allumeur de r\'everb\`eres pour lequel 
$G$ est une HNN extention. Ce r\'esultat affirme que la mesure spectrale
de cet op\'erateur est discr\`ete et concentr\'ee sur un sous ensemble dense de
$[-1,1]$ avec des sauts de valeurs $1/(2^q-1)$, $q \in \mathbb{N}$. 
Les r\'esultats impliquent que $\dim \ker A = \frac{1}{3}$. Mais le d\'enominateur 3
ne divise pas les puissances de 2 qui sont les ordres des sous-groupes finis
du groupe de l'allumeur de r\'everb\`eres.


\begin{center}
--------------------------------------------------------------------------
\end{center}

Atiyah \cite{Atiyah(1976)} introduced for a closed Riemannian manifold
$(M,g)$ with universal covering $\tilde M$ the analytic $L^2$-Betti numbers
$b^p_{(2)}(M,g)$ which measure the size of the space of harmonic
square-integrable $p$-forms on $\tilde M$. Let $k_p(x,y)$ be the (smooth)
integral kernel of the orthogonal projection 
of all
square integrable forms
onto this subspace. On the diagonal, the fiber-wise trace $\tr_x
k_p(x,x)$ is defined and is invariant under deck transformations. It
therefore defines a smooth function on $M$, and Atiyah sets
$b^p_{(2)}(M,g):= \int_M \tr_x k_p(x,x) \;dx$.
By a result of Dodziuk \cite{Dodziuk(1977)}
this does not depend on the metric.

A priori, the $L^2$-Betti numbers are non-negative real
numbers. However, we can express the Euler characteristic $\chi(M)$,
an integer, in terms of the $L^2$-Betti numbers in the usual way:
\begin{equation*}
  \chi(M)=\sum_{p=0}^\infty (-1)^p b^p_{(2)}(M).
\end{equation*}
If $\pi = \pi_1(M)$ is a finite group, then
the $L^2$-Betti numbers can be expressed in terms of ordinary Betti
numbers as follows:
$b^p_{(2)}(M)=\frac{1}{\abs{\pi}} b^p(\tilde M)$.

This note deals with the following conjecture.
Let $\pi$ be a discrete 
group.
 Denote with
  $\fin(\pi)$
   the additive subgroup of $\rationals$ generated by the inverses of
   the orders of the finite subgroups of $\pi$. Note that
   $\fin(\pi)=\integers$ if and only if $\pi$ is torsion free.

\begin{too_optimistic_conjecture}\label{conj:Atiyah}
If $M$ is a closed Riemannian manifold
   with fundamental group $\pi$, then $b^p_{(2)}(M)\in\fin(\pi)$.
   If $\pi$ is torsion free, this specializes to $b^p_{(2)}(M)\in\integers$.
\end{too_optimistic_conjecture}

In \cite{Atiyah(1976)} it is only asked whether the $L^2$-Betti
numbers are always rationals, and integers if the fundamental group is
torsion free. Later, this question was popularized as the Atiyah
conjecture, and also gradually was made precise in the way we
formulate it in Conjecture \ref{conj:Atiyah},
compare \cite[Conjecture 7.1]{Lott-Lueck(1995)}, \cite[Conjecture
9.8]{Linnell(1998)},
 \cite[Conjecture 5.1]{Lueck(1998)},
\cite[Conjecture 5.1]{Reich(1999)},
\cite[Conjecture 13]{Lueck(2000)}, \cite[Definition
1.2]{Dicks-Schick(2000)}, and talks of many mathematicians. 
Conjecture \ref{conj:Atiyah} is also
suggested by \cite[Theorem 8]{Cohen(1979)}, where it is checked for
$\pi$ abelian.

The conjecture is proved in many important cases, starting with
the class $\mathcal{C}$ of Linnell 
\cite{Linnell(1993)}
which includes
extensions of free groups with elementary amenable quotients,
\cite{Schick(1999b)} for residually torsion-free elementary amenable
groups and poly-free groups. In \cite{Dicks-Schick(2000)} it is proved
that the class of groups for which the Atiyah conjecture holds is
closed under HNN-extensions, as long as $\fin(\pi)$ is discrete.
It follows that it holds for all subgroups of
 one-relator
groups, and for all subgroups of right-angled Coxeter groups.
  Dicks-Schick \cite{Dicks-Schick(2000)} also prove that the class of
  all torsion-free groups for which the
  Atiyah conjecture holds is closed under taking extension by groups in
  a certain large class, namely the smallest class which contains all
  the torsion-free, elementary amenable groups, and contains all the
  free groups, and is closed under taking subgroups, extensions,
  directed unions,  amalgamated free products, and HNN-extensions.

 In all the these positive results one has to make one
additional crucial assumption: \emph{there
is a bound on the orders of finite subgroups of $\pi$,
i.e.~$\fin(\pi)$ is a discrete subset of $\reals$}. 

The following
theorem shows that this
additional assumption is essential, and that without it, the 
Strong Atiyah Conjecture
is wrong: 
\begin{theorem}\label{theo:ex}
  Let the group $G$ be given by the presentation
  \begin{equation*}
  G=\generate{a,t,s\mid a^2=1, [t,s]=1, [t^{-1}at,a]=1,
    s^{-1}as=at^{-1}at}.    
  \end{equation*}
We use the notation $[g,h]=g^{-1}h^{-1}gh$
  for the commutator of $g$ and $h$.

The group $G$ is metabelian and in particular elementary amenable.
Every finite subgroup of $G$ is an elementary abelian
  2-group, in particular the order of every finite subgroup of $G$
  is a power of 2.
There exists a closed Riemannian manifold $(M,g)$ of dimension 7
   with $\pi_1(M)=G$
  such that the third $L^2$-Betti number 
  $$b^3_{(2)}(M,g)= \frac{1}{3}.$$ 

\end{theorem}

  Conjecture \ref{conj:Atiyah} predicts that the denominator is
  a power of $2$ and thus the group $G$ is a counterexample to the
  Strong Atiyah Conjecture.

Observe however, that still there is no example of an irrational $L^2$-Betti
number, and that the group $\fin(G)$ is not discrete.

To prove Theorem \ref{theo:ex} we first study the structure of the
group $G$. Let $H$ denote the lamplighter group $(\oplus_{i\in\integers}
\integers/2)\semiProd\integers$, where the generator of $\integers$
acts on $\oplus_{i\in\integers}\integers/2$ by translation. 
 $H$
is generated by $t\in\integers$ and by $a=(\dots,0,1,0,\dots)\in
\oplus_{i\in\integers} \integers/2$ and has the presentation
\begin{equation*}
H=\generate{a,t\mid a^2=1, [t^{-k}at^k,t^{-n}at^n]=1 \;\forall
  k,n\in\integers}. 
\end{equation*}

\begin{lemma}\label{lem:compute_G}
  Let $\alpha \colon H\to H$ be given by $\alpha(t)=t$ and $\alpha(a)=a
  t^{-1}at$. This defines an injective group homomorphism, and $G$ is
  the ascending HNN-extension of $H$ along $\alpha$.
Moreover $G'$ is isomorphic to a countable direct sum of copies of $\integers/2$.
\end{lemma}
\begin{proof}
 The first assertion can be easily checked. The second part follows
 from the computation in \cite{Baumslag(1972)}. For completeness sake
 we give the argument here:

 Let $V$ be the HNN-extension of $H$ along $\alpha$. Then $V$ has the
 presentation
 \begin{equation*}
   V=\generate{a,t,s \mid a^2=1, [s,t]=1, s^{-1}as=a
     t^{-1}at=[a,t], [t^{-k}at^k ,t^{-n}at^n]=1 \;\forall k, n\in\integers}.
 \end{equation*}
 Obviously, we have a epimorphism $G\onto V$ mapping $a$ to $a$, $s$
 to $s$, and $t$ to $t$. It only remains to show that every relation
 in the given presentation of $V$ follows from the relations of
 $G$. Observe first in $G$ that by conjugation with $t^{-n}$,
 $[t^{-k+n}at^{-n+k},a]=1$
 implies $[t^{-k}at^k,t^{-n}at^{n}]=1$. Moreover, commutativity is
 commutative, i.e.~$[t^{-k}at^k,t^{-n}at^{n}]=1$ implies
 $[t^{-n}at^n,t^{-k}at^k]=1$.  Hence, it remains to prove
 $[t^{-n}at^n,a]=1$ in $G$
 for $n>1$. We will do this by induction on
 $n$. Assume therefore $t^{-j}at^j$ commutes with $t^{-l}at^l$ for
 $0\le j\le l<n$. Conjugate the relation $[t^{-(n-1)}at^{n-1},a]=1$
 with $s$. We obtain 
 \begin{equation}\label{commute}
 1=   [t^{-(n-1)}at^{-1}at t^{n-1} ,at^{-1}at] =[(t^{-(n-1)}at^{n-1})(
 t^{-n}a t^n), a (t^{-1}at)].
\end{equation}
Now observe that by induction $a$ commutes with $a_1:=t^{-1}at$ and with
$a_{n-1}t^{-(n-1)}at^{n-1}$. This second relation also implies (by conjugation
  with $t^{-1}$) that moreover, $t^{-1}at$ commutes with $a_n:=t^{-n}a
  t^n$. Therefore, we can simplify the commutator in \eqref{commute}
  to the desired
  \begin{equation*}
    1=(a_n^{-1}a_{n-1}^{-1})(a_1^{-1}a^{-1})(a_{n-1}a_{n})(aa_{1})=
    a_{n}^{-1}(a_{n-1}^{-1}a_{n-1})(a_1^{-1}a_1)a^{-1}a_n a = [t^{-n}at^n,a].
  \end{equation*}
  By induction we therefore see that $V=G$.

  Using the presentation, we next check that the abelianization of $V$
  is isomorphic to $\integers\times\integers$, and $s,t$ are mapped to
  two free generators, whereas $a$ is mapped to zero. Therefore, $G'$
  is equal to the
  normal subgroup generated by $a$, which is generated by
  $s^{-l}t^{-k}a t^k s^l$, $k,l\in\integers$, $l<0$. All these
  elements are of order $2$, and by conjugation with sufficiently
  high powers of $s$ we see that they all commute. Therefore, $G'$ is
  a vector space over $\integers/2$ with countably many generators,
  and therefore isomorphic to a countable direct sum of copies of
  $\integers/2$. Observe, however, that $G'$ is quite different
  from the base of the HNN-extension $H$. The element $sas^{1}$ is a
  typical example which is not contained in $H$ but in $G'$.
\end{proof}

Since,  by Lemma \ref{lem:compute_G}, $G$ is a two-step HNN-extension of
$\oplus_{i\in\integers}\integers/2$, it follows immediately that all finite subgroups of $G$
are elementary abelian $2$-groups. To prove Theorem \ref{theo:ex}, 
we need to construct $M$.

Let $A \in \mathbb {C}G$.  Then left multiplication by $A$ on
$\mathbb {C}G$ induces a bounded linear operator on the Hilbert
space $l^2(G)$; we shall also let $A$ indicate this operator.
Let $\pr_G \colon l^2(G) \to l^2(G)$ denote the orthogonal
projection onto $\ker (A)$ and let $e$ denote the identity element
of $G$; we shall also consider $e$ as the corresponding element of
$l^2(G)$.  We now define $\dim_G(\ker A)$
according to the formula
\[
\dim_G(\ker A) = \langle \pr_G(e), e \rangle _{l^2(G)}.
\]

We will crucially use the following result from \cite[Theorem 2 and Corollary
3]{Grigorchuk-Zuk(2000)}:
\begin{theorem}\label{GriZuk}
  Let $A:= t+at+t^{-1}+(at)^{-1}\in \integers H$ be a multiple of the
Markov
  operator of $H$. Then $A$, considered as an operator on
  $l^2(H)$, has eigenvalues 
  \begin{equation*}
\{ 4\cos\left(\frac{p}{q}\pi\right)\mid
  p\in\integers, q=2,3,\dots\}.
\end{equation*}
The $L^2$-dimension of the
  corresponding eigenspaces is
  \begin{equation*}
    \dim_H \ker\left(A-4\cos\left(\frac{p}{q}\pi\right)\right) =
    \frac{1}{2^q-1}\qquad \text{if }p,q\in\mathbb {Z},
\ q\ge 2, \text{ with }(p,q)=1.
  \end{equation*}
\end{theorem}

  Let us sketch the proof of this theorem.
 In \cite{Grigorchuk-Zuk(2000)} the group $H$ 
 is realized as a group
  defined by a two-state automaton
  (a general
  survey about automata and their groups can be found in
  \cite{Grigorchuk-Nekrashevych-Sushchansky(2000)}).
   Correspondingly, $H$ acts on a
  binary tree and on the boundary of this tree. As to be expected for
  automata, this action shows a lot of self-similarity.
  This
  fractalness can be used to inductively compute the spectra (with
  multiplicity of the eigenvalues) of finite dimensional
  approximations $A_n$ of the operator $A$. The $A_n$ are obtained by
  restricting the action of $H$ to a finite subtree consisting of vertices
  up to the level $n$.  It is important that there exists an infinite
  path in this tree with a trivial stabilizer.
  Similar
  computations have been carried out in
  \cite{Bartholdi-Grigorchuk(1999)}. 
  Then by
  \cite{Grigorchuk-Zuk(1999)}  
  and using approximation
  results for $L^2$-Betti numbers \cite{Farber(1998), Lueck(1994)},
  the spectra (and multiplicities) of the $A_n$ converge, suitably
  normalized, to the $L^2$-dimensions of eigenspaces which have to be
  computed. Observe that in 
  \cite[Corollary 3]{Grigorchuk-Zuk(2000)}, the jump of the spectral
  measure at $4\cos(\frac{p}{q}\pi)$ is computed to be
  $\frac{1}{2^q-1}$, if $(p,q)=1$. Since this jump is exactly the
  $L^2$-dimension of the corresponding eigenspace, the sketch of the
  proof of the
  Theorem is finished.\hfill\qed

  \begin{remark}
    An independent proof of the result of Grigorchuk-Zuk \ref{GriZuk},
    which does not use automata, is given in \cite{Schick(2000b)}.
  \end{remark}

As a corollary of Theorem \ref{GriZuk}, we obtain:
\begin{corollary}\label{corol:matrix}
  There is an $A\in\integers G$ such that 
  \begin{equation*}
\dim_G
  \ker\left(A \colon l^2(G)\to l^2(G)\right) = \frac{1}{3}.    
\end{equation*}
\end{corollary}
\begin{proof}
Observe that if $A$ is induced from $H$, i.e.~$A \in \mathbb{C}H$ (so
we can view $A$ also as an
operator on $l^2(H)$), then essentially $\pr_H  = \pr_G$ and we
deduce that $\dim_H \ker(A) = \dim_G \ker(A)$
(cf.~\cite[Proposition 3.1]{Schick(1998a)}).
Therefore, it will be sufficient to find $A
\in \mathbb {Z}H$ such that $\dim_H \ker A = 1/3$.

Take $A$ of Theorem \ref{GriZuk}. Choosing $p=1$ and $q=2$, we  see
that 0 is in the spectrum of $A$, and that $\dim_H(\ker A) = 1/3$.
\end{proof}

\begin{proposition}\label{prop:CW}
  There is a $3$-dimensional finite CW-complex $X$ with $\pi_1(X)=G$
  and with $b_3^{(2)}(X)=\frac{1}{3}$.
\end{proposition}
\begin{proof}
  We perform a standard construction where one attaching map
  will be given by the $A$ of Corollary \ref{corol:matrix}, compare
  e.g.~\cite[Lemma 2.2]{Lueck(1997)}.

  Let $X'$ be a finite $2$-dimensional CW-complex with $\pi_1(X')=G$,
  e.g. the $2$-complex of the finite presentation given above. Let
  $X''$ be the wedge product of $X'$ and $S^2$. The corresponding
  map $\alpha \colon S^2\to X''$ generates a free copy of
  $\integers\pi_1(X'')=\integers G$ inside $\pi_2(X'')$. Define now
  $X:= X''\cup_{f}D^3$, where $(f \colon S^2\to X'')\in \pi_2(X'')$ is given
  by $A\in\integers G$ of Corollary \ref{corol:matrix}, and where $\integers
  G\into \pi_2(X'')$ is given using $\alpha$. Choosing an appropriate
  basis of cells, it follows that on the cellular $L^2$-chain complex
  $C_*^{(2)}(\tilde X)=C_*(\tilde X)\tensor_{\integers G}l^2(G)$ of
  the universal covering $\tilde X$ of $X$, the differential $d_3$
  \begin{equation*}
    l^2(G) \iso C_3^{(2)}(\tilde X)\overset{d_3}{\longrightarrow}
C_2^{(2)}(\tilde
    X)\iso (l^2(G))^n
  \end{equation*}
  is given by the matrix $(A,0,\dots,0)^t$, where $t$ denotes
transpose and $n$ is the number of $2$-cells
in $X$. Since there are no $4$-cells, $d_4$ is zero. Consequently, 
  \begin{equation*}
    b_3^{(2)}(X) = \dim_G(\ker d_3) = \dim_G(\ker A) = \frac{1}{3}.\qquad\qquad\qed
  \end{equation*}
\renewcommand{\qed}{}
\end{proof}

We now can finish the proof of Theorem \ref{theo:ex} in a standard way
(compare e.g.~\cite[Lemma 2.2]{Lueck(1997)}):

\begin{theorem}\label{theo:manifold}
  There is a $7$-dimensional smooth Riemannian manifold $(M,g)$ with 
  \begin{equation*}
b_3^{(2)}(M,g)=b_3^{(2)}(M)=\frac{1}{3}    
\end{equation*}
and with $\pi_1(M)=G$. Here, $b_3^{(2)}(M)$ denotes the combinatorial
$L^2$-Betti number of a triangulation of $M$.
\end{theorem}
\begin{proof}
  Choose a finite $3$-dimensional simplicial complex $Y$ homotopy
  equivalent to the CW-complex $X$
of Proposition \ref{prop:CW}.
Then embed $Y$ into $\reals^8$ 
\cite[Theorem 5]{Pontryagin(1952)} 
and thicken $Y$ to a homotopy equivalent
  $8$-dimensional compact smooth manifold $W$ with boundary $M$, such
  that moreover the inclusion of $M$ into $W-Y$ is a homotopy equivalence
\cite[Chapter 3]{Rourke-Sanderson(1982)}. Recall that a map $f:V\to V'$
between two CW-complexes is called an $r$-equivalence
($r\in\naturals$), if $\pi_j(f):\pi_j(V)\to \pi_j(V')$ is an
isomorphism for $0\le j< r$, and an epimorphism for $j=r$.
 By transversality (\cite[5.3 and 5.4]{Rourke-Sanderson(1982)}), every
 map of $S^j$ or $D^j$ to $W$ with $j\le 4$ is homotopic to a map into
 $W-Y$. It follows that the
  inclusion $M\hookrightarrow W$ is a $4$-equivalence.  Consequently, by
  \cite[Theorem 1.7]{Lueck(1997)},
  $b_3^{(2)}(M)=b_3^{(2)}(W)=b_3^{(2)}(X)=\frac{1}{3}$
 and
  $\pi_1(M)=\pi_1(X)=G$. If we choose a smooth Riemannian metric $g$
  on $M$, then by the $L^2$-Hodge de Rham theorem \cite[Theorem 1]
{Dodziuk(1977)} we also obtain $b_3^{(2)}(M,g)=\frac{1}{3}$.
\end{proof}

\begin{remark}
  The dimension of the manifold which is a counterexample to the
  strong Atiyah conjecture can be reduced to $6$ as follows:

By \cite[Theorem 3.3]{Lueck(1998b)} $\bigoplus_{i\in\integers}\integers/2$,
$H$, $G$, and the direct limit of $H\xrightarrow{\alpha}
H\xrightarrow{\alpha} \cdots$ all have vanishing $L^2$-Betti
numbers in all degrees. Moreover, the zero-th and first $L^2$-Betti number of a
space is equal to the first $L^2$-Betti number of its fundamental
group \cite[Theorem 1.7]{Lueck(1997)},
i.e.~$b_1^{(2)}(G)=b_1^{(2)}(X)=b_1^{(2)}(M)=0=b_0^{(2)}(X)=b_0^{(2)}(M)$.

The CW-complex $X$ has $1$ zero-cell, $1$ three-cell, and $3$
one-cells and $5$ two-cells (using the presentation of $G$ given in
Theorem \ref{theo:ex}). Consequently, $\chi(X)=2$. Since 
\begin{equation*}
 2= \chi(X) = \sum_{k=0}^3 (-1)^k b^k_{(2)}(X) =
 b^2_{(2)}(X)-b^3_{(2)}(X) b^2_{(2)}(X)-1/3,
    \end{equation*}
we have $b_{(2)}^2(X)=7/3$.

Now we can do the same construction as in the proof of Theorem
\ref{theo:manifold}, but embed $Y$ into $\reals^7$ instead of
$\reals^8$. The inclusion of the boundary $M'$ of the regular
neighborhood $W'$ into $W'$ will now only be a $3$-equivalence, but this
is enough to conclude that $b^2_{(2)}(M')=b^2_{(2)}(W')=b^2_{(2)}(X)
    =7/3$, and the denominator still is not a power of $2$, giving the
    desired counterexample.

Using the K\"unneth formula and Poincar\'e duality \cite[Theorem 1.7]{Lueck(1997)} for
$L^2$-cohomology, one can on the other hand easily
 arrange that the dimension of a counterexample, as well as the
 degree of the Betti number which contradicts the strong Atiyah conjecture, is
arbitrarily high.
\end{remark}

We conclude this paper with a list of open questions regarding the
Atiyah conjecture:
\begin{itemize}
\item Is there an example of an $L^2$-Betti number of
a closed manifold which is not
  rational?
\item Is there an example of a closed manifold with a fundamental
 group $\pi$ with $\fin(\pi)$ discrete
in $\mathbb{R}$
  which provides a counterexample to Conjecture 1?
\item Is there even a counterexample to the Atiyah conjecture with
  torsion free fundamental group? It is well known that, for a torsion
  free group $\pi$,
  the Atiyah conjecture implies that there are no non-trivial
  zero-divisors in $\rationals[\pi]$, even stronger, that
  $\rationals[\pi]$ embeds into a skew field (compare e.g.~\cite[Lemma
  4.4]{Schick(1999b)}). A torsion-free counterexample hence would be
  particularly interesting in view of this zero-divisor conjecture.
\end{itemize}

In the construction of the present counterexample an important role
was played by the
spectral properties of a Markov operator. We would like to formulate the 
following 
open problems:
\begin{itemize}
\item Is there a torsion free group with a Markov operator which 
has a gap in its
spectrum?

\item Is there a torsion free group with a Markov operator whose
 spectral measure is not absolutely continuous with respect to the
 Lebesgue measure or even is not a continuous measure?

\end{itemize} 
These questions are also interesting for the operators given
by any self-adjoint element in a group ring.
{\footnotesize
\bibliographystyle{alpha}
\bibliography{long_counterex}

Preprints of SFB 478, M\"unster are available via http://wwwmath.uni-muenster.de/sfb/about/publ/index.html
}

\end{document}